\documentclass[twoside,reqno]{amsart}

\usepackage{latexsym}

\newcommand{\qdn}{\hspace*{-1.5mm}}
\newcommand{\qqdn}{\hspace*{-2.5mm}}
\newcommand{\xqdn}{\hspace*{-5.0mm}}
\newcommand{\xxqdn}{\hspace*{-10mm}}

\newcommand{\fns}{\footnotesize}
\newcommand{\sst}{\scriptstyle}

\newcommand{\fnk}[3]{\left[\qdn\ba{#1}#2\\#3\ea\qdn\right]}
\newcommand{\ffnk}[4]{\left[\qdn\ba{#1}#3\\#4\ea{\!\Big|\:#2}\right]}

\newcommand{\nnm}{\nonumber}
\newcommand{\be}{\begin{equation}}
\newcommand{\ee}{\end{equation}}
\newcommand{\ba}{\begin{array}}
\newcommand{\ea}{\end{array}}
\newcommand{\bmn}{\begin{eqnarray}}
\newcommand{\emn}{\end{eqnarray}}
\newcommand{\bnm}{\begin{eqnarray*}}
\newcommand{\enm}{\end{eqnarray*}}
\newcommand{\bln}{\begin{subequations}}
\newcommand{\eln}{\end{subequations}}

\newtheorem{thm}{Theorem}

\newtheorem{entry}{Entry}

\newcommand{\bbtm}[4]{\bibitem{kn:#1}{#2,}~\emph{#3,}~{#4.}}
\newcommand{\cito}[1]{\cite{kn:#1}}
\newcommand{\citu}[2]{\cite[#2]{kn:#1}}

\begin{document} 
{\fns
\title{A short proof for Dougall's $_2H_2$-series identity}
\author{$^a$Chuanan Wei, $^b$Qinglun Yan, $^c$ Jianbo Li}
\dedicatory{
$^A$Department of Information Technology\\
  Hainan Medical College,  Haikou 571101, China\\
         $^B$College of Mathematics and Physics\\
   Nanjing University of Posts and Telecommunications,
    Nanjing 210046, China\\
  $^C$Institute of Mathematical Sciences\\
 Xuzhou Normal University, Xuzhou 221116, China}
\thanks{\emph{Corresponding author}: Chuanan Wei. \emph{Email address}:
      weichuanan@yahoo.com.cn}

\address{ }

\begin{abstract}
According to Chen and Fu's method, we offer a short proof for
Dougall's $_2H_2$-series identity.
\end{abstract}

\maketitle\thispagestyle{empty}
\markboth{Chuanan Wei, Qinglun Yan, Jianbo Li}
         {A short proof for Dougall's $_2H_2$-series identity}

Recently, Chen and Fu \cito{chen} gave the semi-finite forms of
several $q$-series identities in a surprising method. Inspired by
this work, we shall derive directly Dougall's $_2H_2$-series
identity form a nonterminating form of Saalsch\"{u}tz's theorem in
the same method.

For a complex number $x$ and an integer $n$, define the shifted
factorial by
 \[(x)_n=\Gamma(x+n)/\Gamma(x)\]
where $\Gamma$-function is well-defined
\[\Gamma(x)=\int_{0}^{\infty}t^{x-1}e^{-t}dt\:\:\text{with}\:\:\mathfrak{Re}(x)>0.\]
For simplifying the expressions, we introduce the following fraction
forms of them:
 \bnm
&&\xqdn\Gamma\fnk{cccc}{a,
&b,&\cdots,&c}{\alpha,&\beta,&\cdots,&\gamma}\:\:\,=\frac{\Gamma(a)\Gamma(b)\cdots
\Gamma(c)}{\Gamma(\alpha)\Gamma(\beta)\cdots \Gamma(\gamma)},\\
&&\xqdn\quad\!\fnk{cccc}{a,&b,&\cdots,&c}{\alpha,&\beta,&\cdots,&\gamma}_n=\frac{(a)_n(b)_n\cdots
(c)_n}{(\alpha)_n(\beta)_n\cdots (\gamma)_n}.
  \enm
 Then a bilateral summation formula due to Dougall
\cito{dougall}(see also Slater \citu{slater}{p. 181}) can be stated
as follows.

\begin{thm}[Dougall's $_2H_2$-series identity]\label{thm-b}
\[\qqdn{_2H_2}\ffnk{ccc}{1}{a,\:b}{c,\:d}=\sum_{k=-\infty}^{\infty}\!\!\fnk{cccc}{a,\:b}
{c,\:d}_k=\Gamma\fnk{cccc}{1-a,\:1-b,\:c,\:d,\:c+d-a-b-1}
{c-a,\:c-b,\:d-a,\:d-b}\] where $\mathfrak{Re}(c+d-a-b)>1$.
\end{thm}

The original proof for Theorem \ref{thm-b} given by Dougall
\cito{dougall} depends on Cauchy residue theorem. For other three
different proofs, the reader may refer to Andrews et al.
\citu{andrews}{p. 110}, Chu \cito{chu} and Slater \citu{slater}{p.
181} respectively. Now, a short new proof for Theorem \ref{thm-b}
will subsequently be offered.

\begin{proof}
Following Bailey~\cito{bailey}, define the unilateral hypergeometric
series by
 \bnm
&&_{1+r}F_s\ffnk{cccc}{z}{a_0,&a_1,&\cdots,&a_r}{&b_1,&\cdots,&b_s}
 \:=\:\sum_{k=0}^\infty
\fnk{cccc}{a_0,&a_1,&\cdots,&a_r}{1,&b_1,&\cdots,&b_s}_kz^k.
 \enm

Then a nonterminating form of Saalsch\"{u}tz's theorem(cf. Andrews
et al. \citu{andrews}{p. 92}) can be expressed as
 \bnm
&&\xxqdn{_3F_2}\ffnk{ccc}{1}{c+d-a-b-1,\:a,\:b}{c,\:d}=
{_3F_2}\ffnk{ccc}{1}{c-a,\:c-b,\:1}{c-a-b+1,\:c+d-a-b}\\
\\&&\xxqdn\:\:\times\:\:\,\Gamma\fnk{cccc}{ c,\:d}
{a,\:b,\:c+d-a-b}\frac{1}{a+b-c}
 +\Gamma\fnk{cccc}{c,\:d,\:c-a-b,\:d-a-b}
{c-a,\:c-b,\:d-a,\:d-b}
 \enm
where $\mathfrak{Re}(d-a-b)>0$. Performing the replacements $k\to
k+n$, $a\to a-n$, $b\to b-n$, $c\to c-n$ and $d\to d-n$ for the last
equation where $k$ denotes the summation index of the $_3F_2$-series
on the left hand side, we have
 \bnm
&&\xqdn\sum_{k=-n}^{\infty}\fnk{cccc}{c+d-a-b-1,\:a-n,\:b-n}
{1,\:c-n,\:d-n}_{k+n}={_3F_2}\ffnk{ccc}{1}{c-a,\:c-b,\:1}{c-a-b+1+n,\:c+d-a-b}\\
\\&&\xqdn\:\:\times\:\:\,\Gamma\fnk{cccc}{\sst c-n,\:d-n}
{\sst a-n,\:b-n,\:c+d-a-b}\frac{1}{a+b-c-n}
 +\Gamma\fnk{cccc}{\sst c-n,\:d-n,\:c-a-b+n,\:d-a-b+n}
{\sst c-a,\:c-b,\:d-a,\:d-b}
 \enm
which is equivalent to the identity
 \bmn
&&\xqdn\sum_{k=-n}^{\infty}\fnk{cccc}{c+d-a-b-1+n,\:a,\:b}
{1+n,\:c,\:d}_k={_3F_2}\ffnk{ccc}{1}{c-a,\:c-b,\:1}{c-a-b+1+n,\:c+d-a-b}\nnm\\\nnm
&&\xqdn\:\:\times\:\:\, \Gamma\fnk{cccc}{\sst 1+n}{\sst c+d-a-b+n}
\Gamma\fnk{cccc}{\sst c,\:d} {\sst
a,\:b}\frac{c+d-a-b-1+n}{(a+b-c-n)(c+d-a-b-1)}
 \\&&\xqdn\:\:+\:\:\,\Gamma\fnk{cccc}{\sst 1+n,\:c-a-b+n,\:d-a-b+n}
{\sst1-a+n,\:1-b+n,\:c+d-a-b-1+n}\Gamma\fnk{cccc}{\sst1-a,\:1-b,\:c,\:d,\:c+d-a-b-1}
{\sst c-a,\:c-b,\:d-a,\:d-b}.\label{dou-for}
 \emn
Recall the limiting relation on $\Gamma$-function:
\[\lim_{n\to\infty}\frac{\Gamma(n+x)}{\Gamma(n+y)}n^{y-x}=1.\]

Letting $n\to \infty$ for \eqref{dou-for} and noting that the first
term on the right hand side vanishes under the condition
$\mathfrak{Re}(c+d-a-b)>1$, we get directly Theorem \ref{thm-b} to
complete the proof.
\end{proof}



\end{document}